\documentclass[12pt]{article}

\usepackage{amssymb,amsmath}
\usepackage[dvips]{graphics,graphicx}

\usepackage{hyperref}

\newcommand{\qed}{\hbox{\rlap{$\sqcap$}$\sqcup$}}
\newtheorem{Lemma}{Lemma}[section]
\newtheorem{definition}{Definition}[section]
 \newenvironment{Definition}{\begin{definition}\rm}
                            {\end{definition}}
\newtheorem{example}{Example}[section]
  \newenvironment{Example}{\begin{example}\rm}
                          {\end{example}}
\newtheorem{prop}{Proposition}[section]
\newtheorem{cor}{Corollary}[section]
\newtheorem{theo}{Theorem}[section]
\newenvironment{pf}{\smallskip\noindent{\bf Proof.}}
                      {{\hfill\qed}\medskip}

\begin{document}

\title{On the arithmetic product of combinatorial species}

\author{Manuel Maia, UCV \and
Miguel M\'{e}ndez, IVIC and UCV}
 \maketitle

\abstract{We introduce two new binary operations with
combinatorial species; the arithmetic product and the modified
arithmetic product. The arithmetic product gives combinatorial
meaning to the product of Dirichlet series and to the Lambert
series in the context of species. It allows us to introduce the
notion of multiplicative species, a lifting to the combinatorial
level of the classical notion of multiplicative arithmetic
function. Interesting combinatorial constructions are introduced;
cloned assemblies of structures, hyper-cloned trees, enriched
rectangles, etc. Recent research of Cameron, Gewurz and Merola,
about the product action in the context of oligomorphic groups,
motivated the introduction of the modified arithmetic product. By
using the modified arithmetic product we obtain new enumerative
results. We also generalize and simplify some  results of
Canfield, and Pittel, related to the enumerations of tuples of
partitions with restricted meet.}

\section{Introduction}

 Informally, a {\em combinatorial species} $F$ (see \cite{BLL,Jo1})
 is a class of labelled
combinatorial structures that is closed by change of labels. Being
more formal, $F$ is a rule assigning to each finite set $U$, a
finite set $F[U]$. The elements of $F[U]$ are called
$F$-structures on the set $U$. The rule $F$ not only acts on
finite sets but also on bijections between finite sets.  To each
bijection $\sigma:U\longrightarrow V,$ the rule $F$ associates a
bijection $F[\sigma]:F[U]\longrightarrow F[V]$ that is called the
\emph{transport of
 $F$-structures along} $\sigma$. In other words, $F$ is an endofunctor
of the category $\mathbb{B}$ of finite sets and bijections.

 For two species of structures $F$ and $G,$ other species can be constructed
 throughout combinatorial operations; addition $F+G,$ product $F\cdot
 G,$ cartesian product $F\times G,$
substitution $F\circ G$ and derivative $F'.$ See \cite{BLL} for
details.

 To each species $F$ are associated three main series expansions.
 The \emph{exponential generating series},
 \begin{equation}\label{F}
    F(x)=\sum_{n\geq 0}|F[n]|\frac{x^n}{n!},
\end{equation}
where $|F[n]|$ is the number of $F$-structures on the set
$[n]=\{1,2,\dots\}.$ The \emph{isomorphism types generating
series},
\begin{equation}\label{Ftilde}
    \widetilde{F}(x)=\sum_{n\geq 0}|F[n]/\sim|x^n,
\end{equation}
where $F[n]/\sim$ denotes the set of isomorphism types of
$F$-structures on $[n].$ The \emph{cycle index series},
 \begin{equation}\label{zetas}
    Z_F(x_1,x_2,\ldots)=\sum_{n\geq 0}\frac{1}{n!}
    \sum_{\sigma\in S_n}\mathrm{fix}F[\sigma]x_1^{\sigma_1}x_2^{\sigma_2}\cdots,
\end{equation}
Here $S_n$ denotes the symmetric group,
$\mathrm{fix}F[\sigma]:=|\mathrm{Fix}F[\sigma]|,$ where
$\mathrm{Fix}F[\sigma]$ is the set of $F$-structures on $[n]$ left
fixed by the permutation $F[\sigma],$ and
 $\sigma_k$ is the number of cycles of length $k$ of $\sigma.$

 Let $F$ be a species of structures and $n\geq 0$ any integer.
 Unless otherwise be explicitly stated, we will denote
 by $F_n$ the species $F$ concentrated in cardinality $n,$
\begin{equation}\label{concentrated}
 F_n[U] =\begin{cases}
 F[U],& \text{if $|U|=n,$}\\
 \varnothing,& \text{if $|U|\not=n$},
 \end{cases}
\end{equation}
where $U$ is a finite set.
  We will also use the notation $F_+$ for the species of
  nonempty $F$-structures,
  \begin{equation}\label{nonempty}
 F_+[U] =\begin{cases}
 F[U],& \text{if $|U|\geq 1,$}\\
 \varnothing,& \text{if $|U|=0$}.
 \end{cases}
\end{equation}

Yeh \cite{yeh0,yeh} established the relationship between
operations with actions of finite permutation groups and
operations with species (product, substitution, cartesian product
and derivative), via the decomposition of a species as a sum of
molecular species. For example, consider the product of two
species
\begin{equation}\label{sproduct}
M\cdot N[U]:=\sum_{U_1+U_2=U}M[U_1]\times N[U_2].
\end{equation}
 When
$M$ and $N$ are molecular,  there are permutation groups $H\leq
S_m$ and $K\leq S_n,$  such that $M=\frac{X^m}{H},$ and
$N=\frac{X^n}{K}.$ We have that
\begin{equation}
\frac{X^m}{H}\cdot\frac{X^n}{K}=\frac{X^{m+n}}{H\times K}
\end{equation}
 where the direct product $H\times K$
acts naturally over the disjoint union $[m]+[n]\equiv[m+n],$ in
what is called the intransitive action. There is another natural
action $H\times K:[m]\times[n],$ the product action, without a
species counterpart. Some enumerative problems have been solved by
Harary \cite{Ha} and by Harrison and High \cite{Harr} using the
cycle index polynomial of the product action.

We can define the arithmetic product of two molecular species by
the formula
\begin{equation}
\frac{X^m}{H}\boxdot\frac{X^n}{K}=\frac{X^{mn}}{H\times K},
\end{equation}
where the action of $H\times K$ over $[mn]\equiv[m]\times[n]$ is
the product action. Then we can extend this product by linearity.
But, in order to have set theoretical definition for the
arithmetic product like formula (\ref{sproduct}) for the ordinary
product, we need a notion of decomposition of a set into factors.
In other words, a set-theoretical analogous of the factoring of a
positive integer as a product of two positive integers. In this
way we arrived to the concept of {\em rectangle} on a finite set.
This concept was previously introduced in other context with the
name of {\em cartesian decomposition} \cite{BPS}, and is a
particular kind of what is called in  \cite{NR} a {\em small
transversal of a partition}.

The most interesting combinatorial construction associated to the
arithmetic product is the {\em assembly of cloned structures}.
Informally, an assembly of cloned $N$-structures is an assembly of
$N$-structures in the ordinary sense, where in addition, all
structures in the assembly are isomorphic replicas of the same
structure. Moreover, information about `homologous vertices' or
`genetic similarity' between each pair in the assembly is also
provided. The structures of $M\boxdot N$ have some resemblance
with the structures of the substitution $M(N).$ An element of
$M\boxdot N$ can be represented as a cloned assembly of
$N$-structures together with an external $M$-structure (an
$M$-assembly of cloned $N$-structures). Because of the symmetry
$M\boxdot N=N\boxdot M$ it also can be represented as an
$N$-assembly of cloned $M$-structures. For example, for $M$ an
arbitrary species  and $L_+$ the species of non-empty lists, the
structures of $M\boxdot L_+$ could be thought of either as
$M$-assemblies of cloned lists, or as lists of cloned
$M$-structures.

There is a link between {\em oligomorphic groups} \cite{Ca2} and
combinatorial species, implicit in the work of Cameron, and which
we hope to have made explicit here. To each oligomorphic group $G$
we can associate a combinatorial species $F_G.$ There is a
correspondence between operations with oligomorphic groups and
operations with species that is very similar to that established
by Yeh between finite permutation groups and molecular species.
For example, the intransitive product action of two oligomorphic
groups translates to the ordinary product of the respective
species and the wreath product to the operation of substitution.
Recently, Cameron,  Gewurz, and Merola have studied the product
action of oligomorphic groups (see \cite{Ca,Ge,GeMe}). This have
motivated us to introduce the {\em modified arithmetic product} in
order to have the appropriated correspondence between operations.
We have made use of this operation to obtain many new enumerative
results. Using a simple manipulation of generating series (the
{\em shift trick}) we greatly simplified and generalized some
results of Canfield \cite{Can}, and Pittel  \cite{Pi}.

\section{The arithmetic product}
\begin{Definition}\label{def1}
For a finite set $U,$ we say that an ordered pair $(\pi,\tau)$ of
partitions of $U$ is a \emph{partial rectangle} on $U$
 when $\pi\wedge\tau=\hat{0}.$ If moreover $\pi$ and $\tau$ are
 independent partitions (every block of $\pi$ meets every block of $\tau$)
 we call it a \emph{rectangle}. More generally, a partial rectangle of
 dimension $k$, or a $k$-partial rectangle is a tuple $(\pi_1,\pi_2,\dots,\pi_k)$
 of partitions such that $\pi_1\wedge\pi_2\wedge\dots\wedge\pi_k=\hat{0}$. It is
called a $k$-rectangle if,
 \begin{equation}
 |B_1\cap\dots\cap B_k|=1,\;\text{ for all }
 B_1\in\pi_1,\dots,B_k\in \pi_k.
 \end{equation}
 \end{Definition}

\begin{figure}[ht]
  \centering
   \includegraphics[height=5cm]{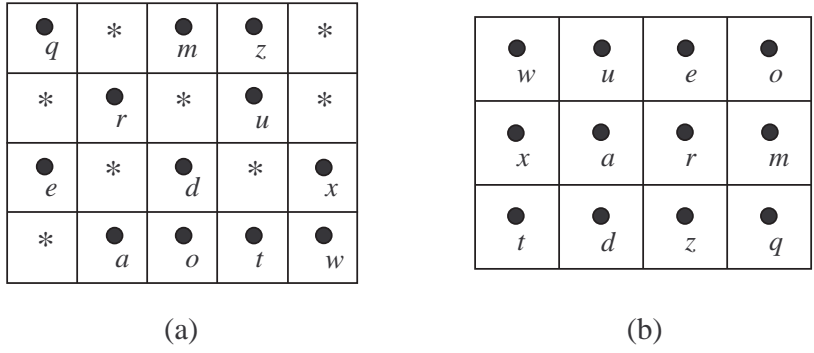}
   \caption{(a) a partial rectangle on a $12$-set, (b) a rectangle on the same set.}
   \label{rectpar}
 \end{figure}

This definition of rectangle is equivalent to the ``cartesian
decomposition" of Baddeley, Praeger and Schneider \cite{BPS}. If
$(\pi,\tau)$ is a rectangle on $U,$ we can arrange the elements of
$U$ in a matrix whose rows are the blocks of $\pi$ and whose
columns are the blocks of $\tau.$ Two matrices represent the same
rectangle if we can obtain one from the other by interchanges of
rows or columns. The same can be say about the partial rectangles
except for the fact that some of the entries of the matrix could
be empty. Figure \ref{rectpar} shows an example of partial
rectangle and rectangle on a set with 12 elements (the symbol
$\ast$ stands by a empty intersection).

 For a rectangle $(\pi,\tau)$ on $U$ obviously holds $|U|=|\pi||\tau|.$
 The \emph{height} of a rectangle $(\pi,\tau)$ is $|\pi|.$ Naturally
$\textrm{height}(\pi,\tau)$ divides $|U|.$ For $|U|=n$ represent
the number of rectangles of  height $d$ with the symbol
$\genfrac{\{}{\}}{0pt}{1}{n}{d}.$ It is no difficult to see that
\begin{equation}
\genfrac{\{}{\}}{0pt}{0}{n}{d}=\displaystyle
\frac{n!}{d!\left(n/d\right)!}.
\end{equation}

Consider $\mathcal{R}$ the species of rectangles, that is, for $U$
a finite set,
  \begin{equation}
   \mathcal{R}[U]=\left\{(\pi,\tau)\mid (\pi,\tau)\mbox{ is a rectangle on } U\right\}.
  \end{equation}
If $n\geq 1,$ we have
  \begin{equation}
   |\mathcal{R}[n]|=\sum_{d|n}\genfrac{\{}{\}}{0pt}{0}{n}{d}.
  \end{equation}
In  an analogous way, for the species $\mathcal{R}^{(k)},$ of
$k$-rectangles, we have
\begin{equation}\label{krectangulo}
   |\mathcal{R}^{(k)}[n]|=\sum_{d_1d_2\cdots
   d_k=n}\genfrac{\{}{\}}{0pt}{0}{n}{d_1,d_2,\dots,d_k},
  \end{equation}
  where
\begin{equation}\label{cor}
  \genfrac{\{}{\}}{0pt}{0}{n}{d_1,d_2,\dots,d_k}=\frac{n!}{d_1!d_2!\cdots
  d_k!}.
\end{equation}

\begin{Definition}\textsc{(Arithmetic product of species)}\label{def2}
 Let $M$ and $N$ be species of structures such that
 $M[\varnothing]=N[\varnothing]=\varnothing.$
 The arithmetic product of $M$ and $N$, is defined as follows
  \begin{equation}
   \left(M\boxdot N\right)[U]:=\sum_{(\pi,\tau)\in\mathcal{R}[U]}M[\pi]\times N[\tau],
  \end{equation}
 where the sum represents the disjoint union and $U$ is a finite set.
  In words, the elements of $\left(M\boxdot N\right)[U]$ are tuples of the form $(\pi,\tau,m,n)$, where
  $m\in M[\pi]$ and $n\in N[\tau]$.
  Recall that given a bijection $\sigma:U\longrightarrow V$ and a partition $\pi$ of $U$, $\sigma$
  induces the partition $\pi'=\sigma(\pi)=\{\sigma(A)\mid A\in \pi\}$ of $V$ and another bijection
  $\sigma^\pi:\pi\longrightarrow \pi'$, sending $A\mapsto\sigma(A),$ for every $A\in
  \pi.$ Similarly for the partition $\tau.$

The transport along a bijection $\sigma:U\longrightarrow V$ is
carried out by
 setting
  \begin{equation}
  \left(M\boxdot N\right)[\sigma]((\pi,\tau,m,n)) = \left(\pi',\tau',M[\sigma^\pi](m),N[\sigma^\tau](n)\right).
  \end{equation}
\end{Definition}

\begin{figure}[ht]
  \centering
   \includegraphics[height=8cm]{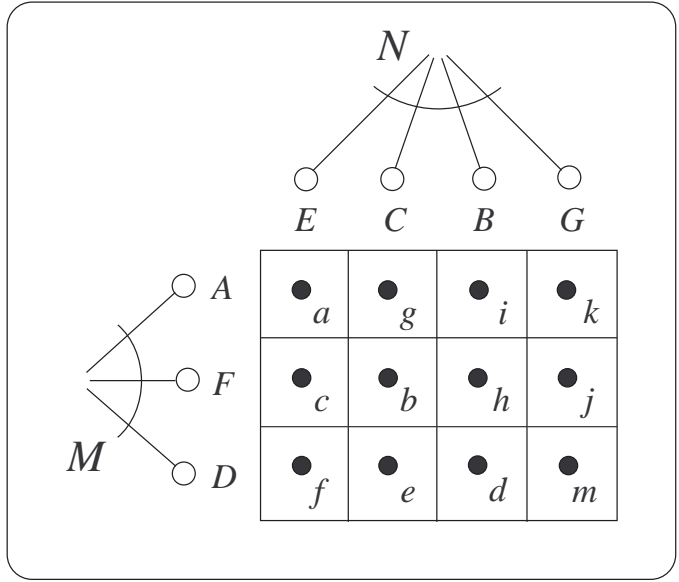}
   \caption{Graphical representation of arithmetic product.}
   \label{prodfig}
 \end{figure}

Figure \ref{prodfig} illustrates an $(M\boxdot N)$-structure on a
set with 12 elements. Here the capital letters (except $M$ and
 $N$) are the labels for the blocks of two partitions forming the rectangle.

\begin{figure}[h]
    \begin{picture}(180,180)
    \put(65,137){\makebox(0,0){$\mathcal{R}$}}
    \put(315,147){\makebox(0,0){$E_+$}}
    \put(258,90){\makebox(0,0){$E_+$}}
    \centering
    \includegraphics[height=6.5cm]{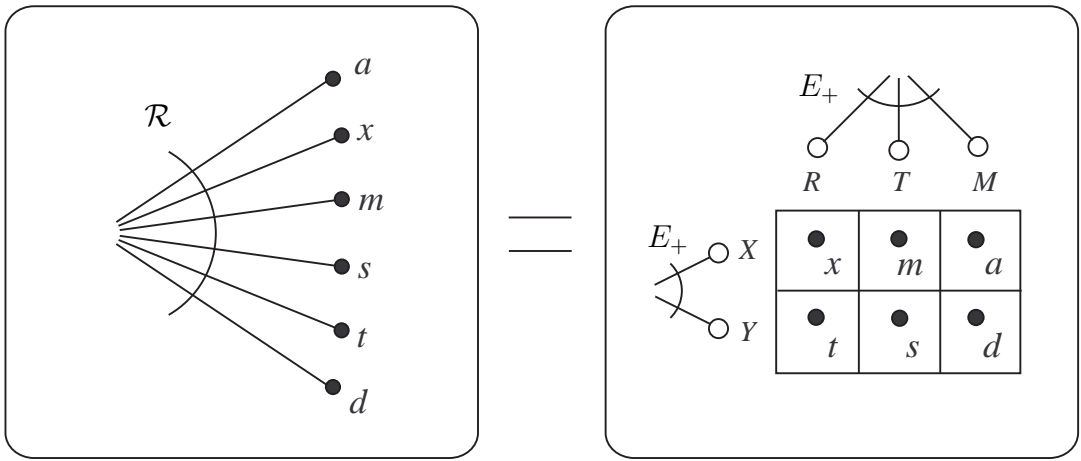}
    \end{picture}
    \caption{An $\mathcal{R}$-structure on a $6$-set.}
    \label{ec_rec}
    \end{figure}

\begin{Example}
Figure \ref{ec_rec} shows that the species $\mathcal{R}$ of
rectangles satisfies the combinatorial equation
$\mathcal{R}=E_+\boxdot E_+,$ where $E_+$ is the species of non
empty sets.
\end{Example}

\begin{prop}\label{card}
Let $M$ and $N$ be species of structures such that
 $M[\varnothing]=N[\varnothing]=\varnothing.$
 Then the exponential generating series of species $M\boxdot N$ is
  \begin{equation}
 (M\boxdot N)(x)=\sum_{n\geq 1}\sum_{d|n}\genfrac{\{}{\}}{0pt}{0}{n}{d}
    |M[d]||N[n/d]|\displaystyle\frac{x^n}{n!}.
\end{equation}
\end{prop}
 \begin{pf}
  \begin{eqnarray*}
   \left|(M\boxdot N)[n]\right| &  =  &\sum_{(\pi,\tau)\in\mathcal{R}[n]}
   \left|M[\pi]\right|\left|N[\tau]\right|\\&=&
   \sum_{d|n}\sum_{\genfrac{}{}{0pt}{1}{(\pi,\tau)\in\mathcal{R}[n]}{\textrm{height}(\pi,\tau)=d}}
   \left|M[d]\right||N[n/d]|\\
    & = & \sum_{d\mid n}\left|\left\{(\pi,\tau)\in\mathcal{R}[n]\mid
   \textrm{height}(\pi,\tau)=d\right\}\right|\left|M[d]\right||N[n/d]|\\
    & = & \sum_{d|n}\genfrac{\{}{\}}{0pt}{0}{n}{d}\left|M[d]\right||N[n/d]|.
  \end{eqnarray*}
\end{pf}

\begin{Example}\textsc{(Regular octopuses}\textrm{ \cite[p. 56]{BLL})}
 Consider the species $\mathcal{C}$ of oriented cycles and
 $L_+$ of non empty linear orders. Figure \ref{rueda}
 represents an
 $(\mathcal{C}\boxdot L_+)$-structure on a set with 8 elements.
 Since $|\mathcal{C}[n]|=(n-1)!$ and $|L_+[n]|=n!,$ we obtain
 \begin{eqnarray*}
   \left|(\mathcal{C}\boxdot L_+)[n]\right| &=&
  \sum_{d|n}\genfrac{\{}{\}}{0pt}{0}{n}{d}|\mathcal{C}[d]||L_+[n/d]|\\
    &=& \sigma(n)(n-1)!,
 \end{eqnarray*}
 where $\sigma(n)$ is the sum of non negative divisors of $n.$
 \begin{figure}[ht]
  \centering
   \includegraphics[height=6.2cm]{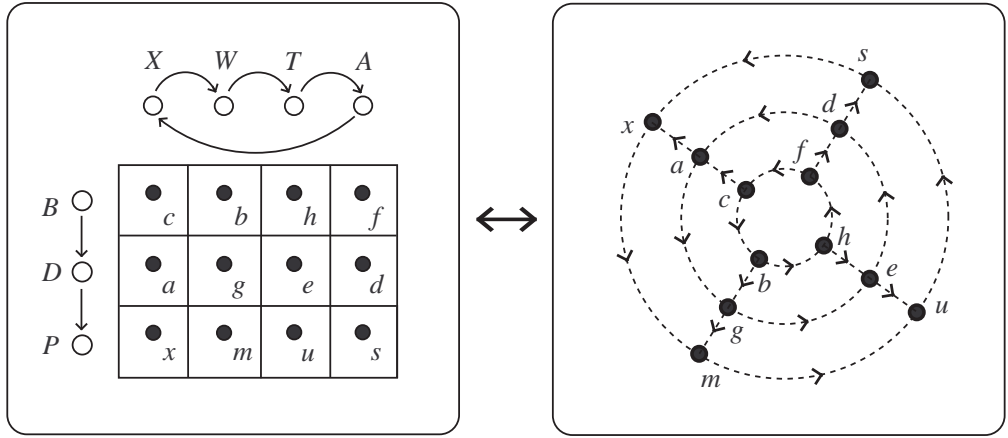}
   \caption{An $(\mathcal{C}\boxdot L_+)$-structure on a $8$-set.}
   \label{rueda}
 \end{figure}
Then, the exponential generating  series is
 \begin{eqnarray*}
   (\mathcal{C}\boxdot L_+)(x) &=& \sum_{n\geq 1}\sigma(n)(n-1)!\frac{x^n}{n!} \\
    &=& \sum_{n\geq 1}\sigma(n)\frac{x^n}{n}.
 \end{eqnarray*}
\end{Example}

\begin{Example}\textsc{(Ordered lists of equal size)}
  Figure \ref{lists} illustrates an
 $(L_+\boxdot L_+)$-structure on a set with 6 elements.
 Since $|L_+[n]|=n!,$ the number of $(L_+\boxdot L_+)$-structures on
 a set of $n$ elements is
 \begin{eqnarray*}
   \left|(L_+\boxdot L_+)[n]\right| &=&
  \sum_{d|n}\genfrac{\{}{\}}{0pt}{0}{n}{d}|L_+[d]||L_+[n/d]|\\
    &=& d(n)n!,
 \end{eqnarray*}
 where $d(n)$ is the number of non negative divisors of $n.$
\begin{figure}[ht]
  \centering
   \includegraphics[height=5.5cm]{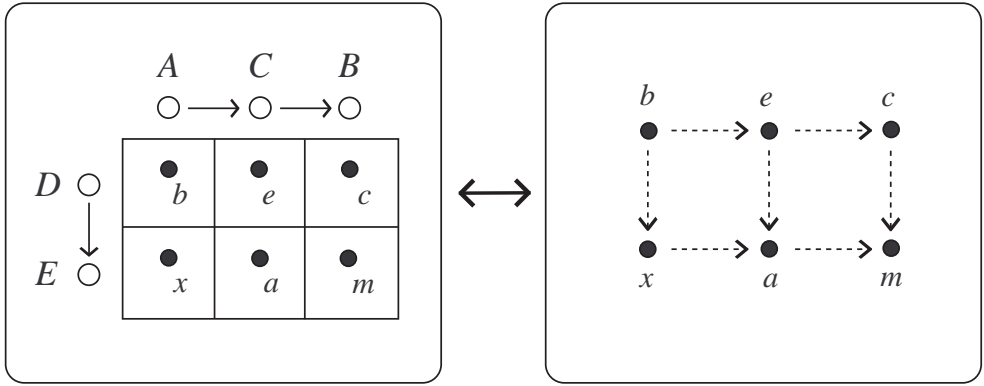}
   \caption{An $(L_+\boxdot L_+)$-structure on a $6$-set.}
   \label{lists}
 \end{figure}
Then, we obtain the generating series
 \begin{eqnarray*}
   (L_+\boxdot L_+)(x) &=& \sum_{n\geq 1} n!d(n)\frac{x^n}{n!} \\
    &=& \sum_{n\geq 1} d(n)x^n.
 \end{eqnarray*}
\end{Example}

\begin{figure}[ht]
  \centering
   \includegraphics[height=7cm]{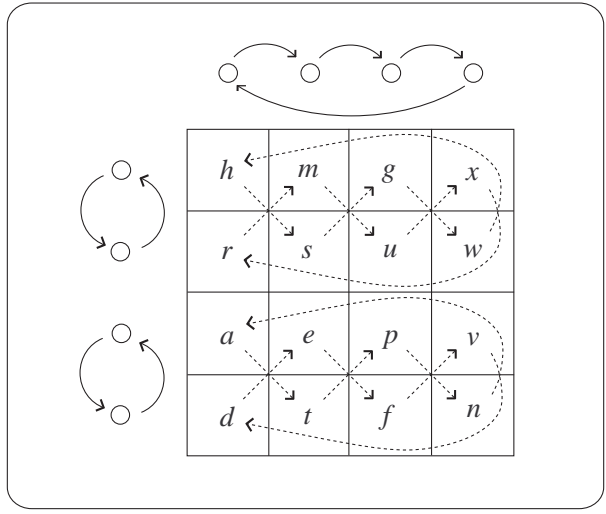}
   \caption{A structure of $\mathcal{S_+}\boxdot \mathcal{S_+}$ and
   induced permutation.}
   \label{biyec}
 \end{figure}
 \begin{Example}\label{marginal}Let $\mathcal{S}_+$ be the species of nonempty permutations.
 It is clear that
 $$
   (\mathcal{S}_+\boxdot \mathcal{S}_+)(x)=(L_+\boxdot L_+)(x)=\sum_{n\geq 1} d(n)x^n.
 $$
The structures of $(\mathcal{S_+}\boxdot \mathcal{S_+})[U]$ are
rectangles enriched with permutations on each side. Formally, they
are tuples of the form $(\pi,\tau,\sigma_1,\sigma_2)$, where
$(\pi,\tau)\in \mathcal{R}[U]$, $\sigma_1\in \mathcal{S_+}[\pi]$,
and $\sigma_2\in \mathcal{S_+}[\tau]$. By the definition of
rectangle, for each element $b\in U$ there exists a unique pair of
sets $(A_b,B_b)\in \pi\times \tau$ such that $b\in A_b\cap B_b$.
The pair $(\sigma _1,\sigma_2)$ induces the permutation
$\sigma_1\boxtimes\sigma_2\in \mathcal{S_+}[U],$ which sends the
element $b\in A_b\cap B_b$ to the unique element in $\sigma
_1(A_b)\cap \sigma_2(B_b)$ (see Figure \ref{biyec}). Let
$\widetilde{\mathcal{R}}$ be the species defined as follows
 \begin{equation}\widetilde{\mathcal{R}}[U]=\{(\pi,\tau,\sigma)\mid
 (\pi,\tau)\in\mathrm{Fix}\mathcal{R}[\sigma], \sigma\in
 \mathcal{S}[U]\}.
 \end{equation}
The function
\begin{eqnarray*}\label{isa}
 \boxtimes_U:(\mathcal{S_+}\boxdot
 \mathcal{S_+})[U]&\longrightarrow&
 \widetilde{\mathcal{R}}[U]\\
 (\pi,\tau,\sigma_1,\sigma_2)&\longmapsto&(\pi,\tau,\sigma_1\boxtimes\sigma_2)
\end{eqnarray*}
is a natural bijection with inverse
$(\pi,\tau,\sigma)\mapsto(\pi,\tau,\sigma^{\pi},\sigma^{\tau})$.
The family $\{\boxtimes_U\}_{U\in \mathbb{B}},$ defines
 a species isomorphism
 $\boxtimes:\mathcal{S_+}\boxdot \mathcal{S_+}\longrightarrow \widetilde{\mathcal{R}}.$
\end{Example}
\begin{prop}\label{prop}
 Let $M,N$ and $R$ be species of structures such that $M[\varnothing]=N[\varnothing]=R[\varnothing]=\varnothing,$
  and $X$ the singular species.
 The product $\boxdot$ has the following properties:
 \begin{eqnarray}
    M\boxdot N&=&N\boxdot M,\\
    \label{asociatividad}M\boxdot(N\boxdot R)&=&(M\boxdot N)\boxdot R,\\
   \label{distributividad} M\boxdot(N+R)&=&M\boxdot N+M\boxdot R,\\
    M\boxdot X&=&X\boxdot M=M,\\
    (M\boxdot N)^{\bullet}&=&M^{\bullet}\boxdot N^{\bullet},\\
    \label{listas}M\boxdot X^n&=&M(X^n),\\
    M\boxdot L_+&=&\sum_{n\geq 1}M(X^n).\label{assemb}
 \end{eqnarray}
 \qed
\end{prop}
All the properties are not difficult to prove. In particular, the
reader may verify that both sides of equation
(\ref{asociatividad}) evaluated at a set $U$, are naturally
equivalent to the set
\begin{equation}
\sum_{(\pi_1,\pi_2,\pi_3)\in \mathcal{R}^{(3)}}M[\pi_1]\times
N[\pi_2]\times R[\pi_3].
\end{equation}
In general, for a family $\{ M_i\}_{i=1}^k$ of species with
$M_i[\varnothing]=\varnothing$, we have
\begin{equation}\label{general}
(\boxdot_{i=1}^kM_i)[U]=\sum_{(\pi_1,\pi_2,\dots,\pi_k)\in
\mathcal{R}^{(k)}[U]}\prod_{i=1}^k M_i[\pi_i],
\end{equation}
for every $i=1,2,\dots,k.$

\begin{figure}[ht]
  \centering
   \includegraphics[height=6cm]{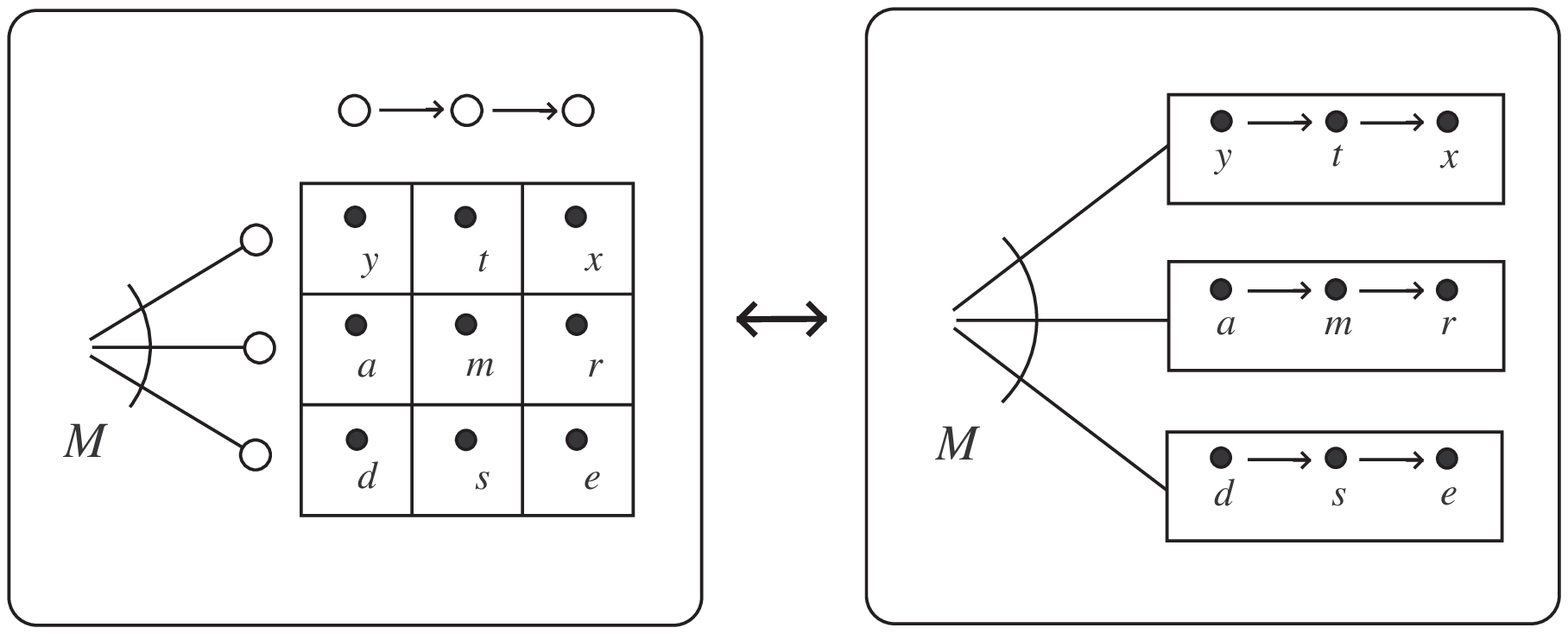}
   \caption{An $(M\boxdot L_+)$-structure.}
   \label{Massemb}
 \end{figure}

 From equation (\ref{assemb}), the
structures of $M\boxdot L_+$ may be thought of as $M$-assemblies
of lists of equal size (see Figure \ref{Massemb}).

\subsection{The arithmetic product and generating series}
\begin{Definition}\label{arithp}
For any two monomials $x^n$ and $x^m$ we define the arithmetic
product $x^n\boxdot x^m:=x^{nm}.$ Extend this product by linearity
to exponential formal power series with zero constant term.
\end{Definition}
We easily obtain that
\begin{equation}
\left(\sum_{n\geq 1}a_n\frac{x^n}{n!}\right)\boxdot
\left(\sum_{n\geq 1}b_n\frac{x^n}{n!}\right)=\sum_{n\geq
1}c_n\frac{x^n}{n!},
\end{equation}
where
\begin{equation}\label{producto aritmetico}
c_n=\sum_{d|n}\genfrac{\{}{\}}{0pt}{0}{n}{d}a_d
b_{n/d}.
\end{equation}

Observe that the exponential formal power series with the
arithmetic product form a ring  with identity $x$. The
substitution $\displaystyle x^n \leftarrow \frac{1}{n^s},$ makes
it isomorphic to the ring of \emph{modified formal Dirichlet
series},
$$\sum_{n\geq 1}\frac{a_n}{n!n^s}.$$
This motivates the following definition.
\begin{Definition}
 Let $M$ be a species of structures satisfying the condition
 $M\left[\varnothing\right]=\varnothing.$ Then the
 \emph{modified Dirichlet generating series} of $M$ is
 \begin{equation}
  \mathcal{D}_M(s)=\sum_{n\geq 1}\frac{|M[n]|}{n!n^{s}}.
 \end{equation}
\end{Definition}
Thus, for the species $E_+,L_+,\mathcal{C}$ and $S_+,$ we have:
\begin{eqnarray}
  \displaystyle\mathcal{D}_{E_+}(s) &=& \sum_{n\geq 1}\frac{1}{n!n^{s}},\\
  \displaystyle\mathcal{D}_{L_+}(s) &=&   \sum_{n\geq 1}\frac{1}{n^{s}}=\zeta(s),\\
  \displaystyle\mathcal{D}_{\mathcal{C}}(s) &=&   \sum_{n\geq 1}\frac{1}{n^{s+1}}=\zeta(s+1),\\
  \displaystyle\mathcal{D}_{S_+}(s) &=&   \sum_{n\geq 1}\frac{1}{n^{s}}=\zeta(s).
\end{eqnarray}

From Proposition \ref{card} we obtain the following.
\begin{prop}\label{prodDir}
For species of structures $M$ and $N$ with the condition
$M\left[\varnothing\right]=\varnothing=N\left[\varnothing\right],$
we have:
\begin{equation}
(M\boxdot N)(x)=M(x)\boxdot N(x),
\end{equation}
and
\begin{equation}
 \mathcal{D}_{M\boxdot
N}(s)=\mathcal{D}_M(s)\cdot \mathcal{D}_N(s).
\end{equation}
\end{prop}

For a formal power series $R(x)$ we have
\begin{equation}
x^n\boxdot
R(x)=R(x^n).
\end{equation}
Thus we have the generating series identity
\begin{equation}\label{noisomorfo}
(M\boxdot N)(x)=\sum_{n\geq 1}\frac{|M[n]|}{n!}N(x^n).
\end{equation}
In particular
\begin{equation}
x^n\boxdot \frac{x}{1-x}=\frac{x^n}{1-x^n},
\end{equation}
and we obtain that the generating series of $M\boxdot L_+$ is the
{\em Lambert series}
\begin{equation}
(M\boxdot L_+)(x)=\sum_{n\geq
1}\frac{|M[n]|}{n!}\frac{x^n}{1-x^n}.
 \end{equation} By (\ref{assemb}) we also have
\begin{equation}
(M\boxdot L_+)(x)=\sum_{n\geq 1}M(x^n)
 \end{equation}
  Using the previous two equations we get
\begin{eqnarray}
(\mathcal{C}\boxdot L_+)(x)&=&\sum_{n\geq
1}\frac{x^n}{n(1-x^n)}=\sum_{n\geq
1}\ln\left(\frac{1}{1-x^n}\right)=\sum_{n\geq
1}\frac{\sigma(n)}{n}x^n,\\
(L_+\boxdot L_+)(x)&=&\sum_{n\geq 1}\frac{x^n}{1-x^n}=\sum_{n\geq
1}d(n)x^n,\\
(L_+^{\bullet}\boxdot L_+)(x)&=&\sum_{n\geq
1}n\frac{x^n}{1-x^n}=\sum_{n\geq
1}\frac{x^n}{(1-x^n)^2}=\sum_{n\geq 1}\sigma(n)x^n.
\end{eqnarray}(see \cite{Com} and \cite{W} for more properties of
Lambert series). Those identities translate to Dirichlet
generating series as
\begin{eqnarray}
\mathcal{D}_{\mathcal{C}\boxdot
L_+}(s)&=&\zeta(s+1)\zeta(s)=\sum_{n\geq
1}\frac{\sigma(n)}{n}n^{-s},\\
\mathcal{D}_{L_+\boxdot L_+}(s)&=&\zeta^2(s) = \sum_{n\geq 1}^{}d(n)n^{-s},\\
\mathcal{D}_{L_+^{\bullet}\boxdot
L_+}(s)&=&\zeta(s-1)\zeta(s)=\sum_{n\geq 1}\sigma(n)n^{-s}.
\end{eqnarray}
By equation (\ref{noisomorfo}) we also obtain:
\begin{eqnarray}\label{unital1}
(\mathcal{C}\boxdot M)(x)=\ln\prod_{n\geq
1}\left(\frac{1}{1-x^n}\right)^{\frac{|M[n]|}{n!}},\\
\label{unital2}E(\mathcal{C}\boxdot M)(x)=\prod_{n\geq
1}\left(\frac{1}{1-x^n}\right)^{\frac{|M[n]|}{n!}}.
\end{eqnarray}

Let $M$ be a species of structures. To describe the compatibility
of the product $\boxdot$ with the transformation $M\rightarrow
Z_M,$ it is necessary to define a product $\boxdot$ for two index
series (see \cite{Ha}). First we have the following Lemma.
\begin{Lemma}Let $(\pi,\tau,\sigma_1,\sigma_2)$ be an element of
$\left(\mathcal{S_+}\boxdot\mathcal{S_+}\right)[U]$ and
$\sigma=\sigma_1\boxtimes\sigma_2\in \mathcal{S}[U]$. If the cycle
type of $\sigma,$ $\sigma_1,$ and $\sigma_2,$  are respectively
$\alpha,$ $\beta$ and $\gamma,$ then we have
\begin{equation}\label{caja}
\alpha_k=\sum_{[i,l]=k}(i,l)\beta_i \gamma_l,\quad
k=1,2,\dots,[d,n/d],
\end{equation}
where $d=\emph{height}(\pi,\tau),$ $[i,l]$ denotes the
least common multiple of\, $i$ and $l,$ and $(i,l)$ the greatest common divisor.
\end{Lemma}
\begin{pf}

Analogous to the proof of proposition 7(b) in (\cite{BLL} p. 74).
\end{pf}

We will say that  $\alpha=\beta\boxtimes\gamma$ when they satisfy
the equation (\ref{caja}). Define now the operation $\boxdot$ on
monomials by
 \begin{equation}
  \left(\prod_{i=1}^n x_i^{\beta_i}\right)\boxdot
   \left(\prod_{l=1}^m x_l^{\gamma_l}\right)
     :=\prod_{i=1}^n\prod_{l=1}^m x_{[i,l]}^{\beta_i\gamma_l(i,l)}
      =\prod_{k=1}^{nm}x_k^{\alpha_k},
 \end{equation}
where $\alpha_k=\sum_{[i,l]=k}(i,l)\beta_i \gamma_l=(\beta
\boxtimes \gamma)_k.$ Equivalently
\begin{equation}\label{productar}
\mathbf{x}^{\beta}\boxdot\mathbf{x}^{\gamma}:=\mathbf{x}^{\beta
\boxtimes \gamma}.
\end{equation}
Finally extend linearly this operation to polynomials and formal
power series. Note that $x_1^i\boxdot x_1^j=x_1^{ij}$ as in
Definition~\ref{arithp}.
\begin{Definition}For $\alpha,$ $\beta,$ and $\gamma$ as above,
define the coefficient
\begin{equation}
\genfrac{\{}{\}}{0pt}{0}{\alpha}{\beta,\gamma}
=
\begin{cases}
\frac{\mathrm{aut}(\alpha)}
{\mathrm{aut}(\beta)\mathrm{aut}(\gamma)},& \text{if }\beta
\boxtimes \gamma=\alpha,\\0,& \text{otherwise}.
\end{cases}
\end{equation}
\end{Definition}
\begin{Lemma}\label{binomial}Let $\sigma$ be a permutation on a finite set $U$ with cycle type
$\alpha$. For $\beta$ and $\gamma$ as above, let
$S_{\beta,\gamma}^{\sigma}$ be the set of tuples
$(\pi,\tau,\sigma_1,\sigma_2)\in
(\mathcal{S}_+\boxdot\mathcal{S}_+)[U],$ such that
$\sigma_1\boxtimes\sigma_2=\sigma$, and the cycle type of
$\sigma_1$ and $\sigma_2$ are respectively $\beta$ and $\gamma.$
Then \begin{equation}
|S_{\beta,\gamma}^{\sigma}|=\genfrac{\{}{\}}{0pt}{0}{\alpha}{\beta,\gamma}.
\end{equation}
\end{Lemma}
\begin{pf}
The group $\mathrm{Aut}(\sigma)$ acts transitively on
$S_{\beta,\gamma}^{\sigma}$ in the following manner: for $\eta\in
\mathrm{Aut}(\sigma),$
\begin{equation}
\eta\cdot(\pi,\tau,\sigma_1,\sigma_2):=(\eta(\pi),\eta(\tau),
\eta^{\pi}\sigma_1(\eta^{\pi})^{-1},
\eta^{\tau}\sigma_2(\eta^{\tau})^{-1}).
\end{equation}
The order of the group fixing any element of
$S_{\beta,\gamma}^{\sigma}$ is
$\mathrm{aut}(\beta)\mathrm{aut}(\gamma).$ \end{pf}
\begin{prop}\label{index}
 Let $M$ and $N$ be two species of structures. Then, the cycle index
 series and the type generating
 series associated to the species $M\boxdot N$ satisfy the identities
 \begin{eqnarray}
    \label{cycle}Z_{M\boxdot N}(x_1,x_2,\ldots)&=&Z_M(x_1,x_2,\ldots)\boxdot Z_N(x_1,x_2,\ldots),\\
  \label{types} \widetilde{M\boxdot N}(x)&=&\widetilde{M}(x)\boxdot \widetilde{N}(x).
 \end{eqnarray}
\end{prop}
\begin{pf}
It is not difficult to deduce the second identity from the first.
Using equation (\ref{productar}) we obtain
\begin{equation}Z_M(\mathbf{x})\boxdot
Z_N(\mathbf{x})=\sum_{\alpha}\left(\sum_{\beta \boxtimes
\gamma=\alpha}\genfrac{\{}{\}}{0pt}{0}{\alpha}{\beta,\gamma}\mathrm{fix}M[\beta]
\mathrm{fix}N[\gamma]\right)\frac{\mathbf{x}^{\alpha}}{\mathrm{aut}(\alpha)}.
\end{equation}
Then, all we have to prove is that
\begin{equation}\mathrm{fix}(M\boxdot N)[\alpha]=|\mathrm{Fix}(M\boxdot
N)[\sigma]|=\sum_{\beta \boxtimes
\gamma=\alpha}\genfrac{\{}{\}}{0pt}{0}{\alpha}{\beta,\gamma}\mathrm{fix}M[\beta]
\mathrm{fix}N[\gamma],
\end{equation}
where $\sigma$ is any permutation on a finite set $U,$ with cycle
type $\alpha.$ Since
\begin{equation}
\mathrm{Fix}(M\boxdot N)[\sigma]=\left\{(\pi,\tau,m,n)\mid
(\pi,\tau)\in\mathrm{Fix}\mathcal{R}[\sigma], m\in
\mathrm{Fix}M[\sigma^{\pi}],
n\in\mathrm{Fix}N[\sigma^{\tau}]\right\},
\end{equation}
we have
\begin{eqnarray}
\mathrm{fix}(M\boxdot
N)[\alpha]&=&\sum_{(\pi,\tau)\in\mathrm{Fix}\mathcal{R}[\sigma]}|\mathrm{Fix}M[\sigma^{\pi}]|
|\mathrm{Fix}N[\sigma^{\tau}]|\\
&=&\sum_{\genfrac{}{}{0pt}{1}{(\pi,\tau,\sigma_1,\sigma_2)\in(\mathcal{S}_+\boxdot\mathcal{S}_+)[U]}
{\sigma_1\boxtimes\sigma_2=\sigma}}|\mathrm{Fix}M[\sigma_1]||\mathrm{Fix}N[\sigma_2]|.
\end{eqnarray}
The last identity is obtained from bijection (\ref{isa}) in
example (\ref{marginal}). Classifying the permutations $\sigma_1$
and $\sigma_2$ according with their cycle type, we get
\begin{equation}
\mathrm{fix}(M\boxdot
N)[\alpha]=\sum_{\beta,\gamma}\sum_{(\pi,\tau,\sigma_1,\sigma_2)\in
S_{\beta,\gamma}^{\sigma}}\mathrm{fix}M[\beta]\mathrm{fix}N[\gamma].
\end{equation}
By lemma (\ref{binomial}) we obtain the result.
\end{pf}

\subsubsection{ The cyclotomic identity}
There are various
bijective proofs of the cyclotomic identity
$$\frac{1}{1-\alpha x}=\prod_{n\geq
1}\left(\frac{1}{1-x^n}\right)^{\lambda_n(\alpha)}.$$ See for
example: Metropolis-Rota \cite{MR}, Taylor \cite{T} and Bergeron
\cite{Ber}. We propose here a very simple one, as an application
of the combinatorics of the arithmetic product.

 Let $\mathcal{C}^{(\alpha)}$ be
the species $\alpha$-colored cycles, or necklaces, following the
terminology of Metropolis and Rota (see \cite{MR}). The elements
of $\mathcal{C}^{(\alpha)}[U]$ are pairs of the form $(\sigma,f)$,
where $\sigma\in \mathcal{C}[U]$ and $f:U\rightarrow \mathrm{A}$
is an arbitrary function assigning colors (letters) in a  totally
ordered set $\mathrm{A}$ (alphabet) with $|\mathrm{A}|=\alpha,$ to
the labelled beads of the cycle $\sigma$. Denote by
$\mathcal{S}^{(\alpha)}$ the species of assemblies of necklaces.
It is clear that
\begin{eqnarray}
\mathcal{C}^{(\alpha)}(x)&=&\ln\left(\frac{1}{1-\alpha x}\right),\\
\mathcal{S}^{(\alpha)}(x)&=&\frac{1}{1-\alpha x}.
\end{eqnarray}
Let $(\sigma,f)$ be a necklace in $\mathcal{C}^{(\alpha)}[U],$
where $|U|=n.$ The integer $$d=\mathrm{min}\{1\leq k\leq n\mid
f\circ\sigma^k=f\}$$ is called the {\em period} of $(\sigma,f). $
When $d=n$, the necklace is called {\em aperiodic}.
 The {\em flat part} $\overline{\mathcal{C}^{(\alpha)}}$ of
$\mathcal{C}^{(\alpha)}$ (see \cite{La1}) is the species of
aperiodic necklaces.

 Let $(\sigma,f)$ be an aperiodic necklace. To each of the $n$
 possible presentations of the
 cycle $\sigma$ as an ordered tuple
 $\sigma=(a_1,a_2,\dots,a_n)$
  corresponds a different word $f(a_1)f(a_2)\dots f(a_n)$ in the alphabet $\rm{A}.$
 The lowest of them in the lexicographic order is  called a
 {\em Lyndon} word. The ordering $\sigma=(a_1,a_2,\dots,a_n)$ such
 that the corresponding word $w$ is Lyndon will be called the {\em
 standard presentation} of $\sigma.$ Thus, the necklace $(\sigma,f)$
 can be identified with the pair $(w,l),$ where $w$ is a Lyndon word
 and $l,$ a linear order on $U$, is the standard presentation of the cycle $\sigma$.
 The number of Lyndon words on
 $\mathrm{A}$ is well known to be
 $\lambda_n(\alpha):=\frac{1}{n}\sum_{d|n}\mu(d)\alpha^{n/d},$
 where $\mu$ is the classical M\"{o}bius function.
 Then, the decomposition of $\overline{\mathcal{C}^{(\alpha)}}$
 as a sum of molecular species is
\begin{equation}
\overline{\mathcal{C}^{(\alpha)}}=\sum_{n\geq
1}\lambda_n(\alpha)X^n.
\end{equation}

\begin{prop}We have the equalities:
\begin{eqnarray} \label{cyclotomic1}
\mathcal{C}^{(\alpha)}&=&\mathcal{C}\boxdot\overline{\mathcal{C}^{(\alpha)}},\\
\label{cyclotomic2}\mathcal{S}^{(\alpha)}&=&E(\mathcal{C}
\boxdot\overline{\mathcal{C}^{(\alpha)}}).
\end{eqnarray}
\end{prop}
\begin{pf}
Equation (\ref{cyclotomic2}) is immediate from
(\ref{cyclotomic1}).
\begin{figure}[!ht]
  \centering
   \includegraphics[height=6.5cm]{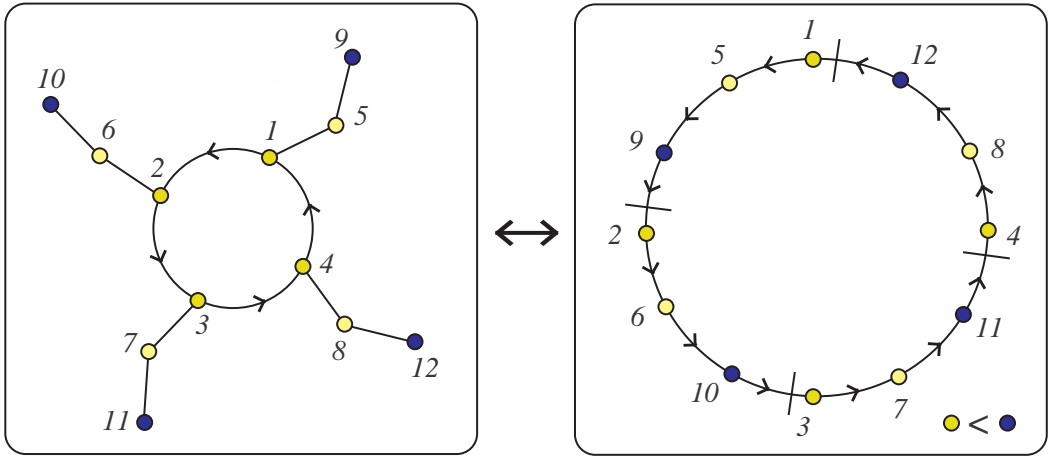}
   \caption{The isomorphism
   $\mathcal{C}^{(\alpha)}=\mathcal{C}\boxdot\overline{\mathcal{C}^{(\alpha)}}.$}
   \label{cyc_fig}
 \end{figure}
To prove (\ref{cyclotomic1}) observe that the structures of
$\mathcal{C}\boxdot \overline{\mathcal{C}^{(\alpha)}}$ are regular
octopuses where each tentacle (linear order) is decorated with the
same Lyndon word. Join the decorated tentacles following the
external cycle of the octopus to obtain a necklace whose period is
the common length of the tentacles. Conversely, given a necklace
of period $d,$ there is a unique way of cutting it into pieces of
length $d$ such that the word on each piece is Lyndon. It is easy
to see how to get an element of
 $\mathcal{C}\boxdot\overline{\mathcal{C}^{(\alpha)}}$ out of this sliced
 necklace. See for example \cite[pages 4-5]{St}, where a similar bijection is
 used to count ordinary octopuses.
\end{pf}

Equation (\ref{cyclotomic2}) can be interpreted as the cyclotomic
identity lifted at a combinatorial level. By equations
(\ref{unital1}) and ({\ref{unital2}), we get:
\begin{eqnarray}
 \mathcal{C}^{(\alpha)}(x)&=&\ln\prod_{n\geq
1}\left(\frac{1}{1-x^n}\right)^{\lambda_n(\alpha)},\\
\mathcal{S}^{(\alpha)}(x)&=&\frac{1}{1-\alpha x}=\prod_{n\geq
1}\left(\frac{1}{1-x^n}\right)^{\lambda_n(\alpha)}.
\end{eqnarray}

\section{Assemblies of cloned structures}
In this section we will see that an $(M\boxdot N)$-structure can
be interpreted as an ``$M$-assembly of cloned $N$-structures". The
intuition behind this is the following: an element of $(M\boxdot
N)[U]$ consist of a rectangle $(\pi,\tau)$ on $U$ enriched with an
$M$-structure $m$ on one side ($\pi$) and an $N$-structure $n$ on
the other side $(\tau).$ Because $\tau$ have the same number of
elements than any block of $\pi,$ we could laid an isomorphic copy
(clone) of $n$ on each block $B$ of $\pi$. Those copies of $n$
together with the ``external structure'' $m\in M[\pi]$ form an
$M$-assembly of cloned $N$-structures.

To make this definition precise we need some formalism. Two
elements of $U$ belonging to the same block of $\tau$ will be
called \emph{homologous}. For example let $\mathcal{A}$ be the
species of rooted trees. In Figure \ref{homol} we represent in two
ways a structure of $\mathcal{C}\boxdot \mathcal{A}$ as a
$\mathcal{C}$-assembly of cloned rooted trees. In both of them
homologous elements are represented with the same color (pattern).
Roots of cloned trees in the right hand side are connected like
the original cycle on $\pi$ in the left hand side, the rest of
homologous elements are connected with closed segmented curves.

 We now express conveniently the relation among
homologous elements.
 Let $B\in\pi$ and $b\in B.$ It is clear that there is only one
block $C\in\tau$ such that $\{b\}=B\cap C.$ For
$\left(B,B'\right)\in\pi\times \pi,$ we define the bijection
 $$
  \begin{array}{cl}
    \Phi_{B,B'}^{\tau} : &\!\!\! B \longrightarrow B'\\
                  &\!\! b\, \longmapsto\, b', \\
  \end{array}
 $$
where $b'$ is the unique element of $B'\cap C.$ In other words,
$\Phi_{B,B'}^{\tau}$ sends each element of $B$ to its homologous
in $B'$. It is easy to verify that
\begin{itemize}
 \item[(i)] $\Phi_{B,B}^{\tau}=\textrm{Id}_B$ and
 \item[(ii)] $\Phi_{B',B''}^{\tau}\circ\Phi_{B,B'}^{\tau}=\Phi_{B,B''}^{\tau},$ for all
 $B,B',B''\in\pi.$
\end{itemize}
\begin{Definition}Let  $U$ be a finite set and $M,N$ two
species of structures such that
$M[\varnothing]=N[\varnothing]=\varnothing.$ An $M$-assembly of
cloned $N$-structures is a triple $(\{n_B\}_{B\in\pi},\tau,m),$
where:
 \begin{itemize}
  \item[(i)] $(\pi,\tau)\in\mathcal{R}[U],$
  \item[(ii)] $\{n_B\}_{B\in\pi}$ is an assembly of
  $N$-structures ($n_B\in N[B],$ for each $B\in\pi$), along with the condition
 \begin{equation}\label{cond}
    N[\Phi_{B,B'}]n_B=n_{B'}, \end{equation}
    for every pair $(B,B')\in\pi\times \pi,$
  \item[(iii)] $m\in M[\pi].$
 \end{itemize}
\end{Definition}
\begin{figure}[ht]
  \centering
   \includegraphics[height=6.7cm]{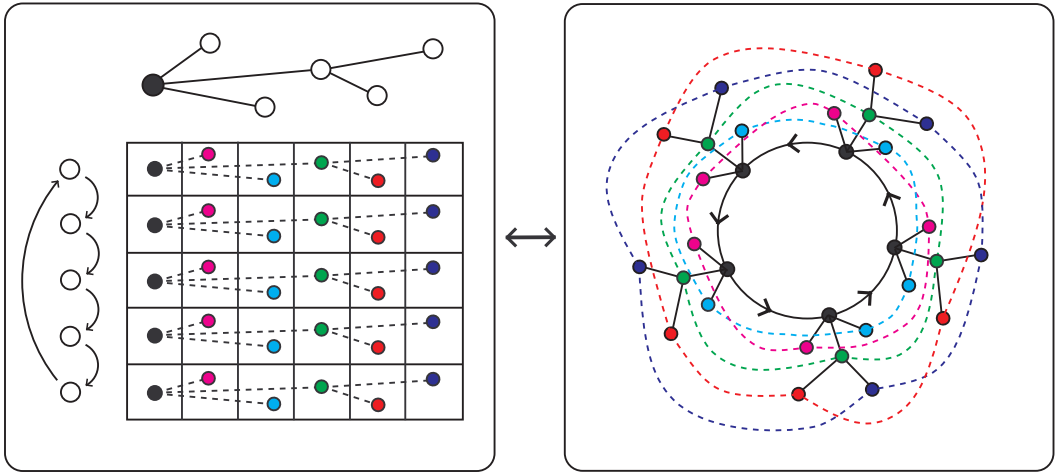}
   \caption{A $\mathcal{C}$-assembly of cloned rooted trees.}
   \label{homol}
 \end{figure}

\begin{prop}\label{cloned}
Let be $M$ and $N$ be two species of structures. Then the species
$M\boxdot N$ and the species of $M$-assemblies of cloned
$N$-structures are isomorphic.
\end{prop}
\begin{pf}
 Let $U$ be a finite set and assume that
 $(\pi,\tau,m,n)\in(M\boxdot N)[U].$
 For each $B\in\pi,$ let $\Psi_{\tau,B}:\tau\longrightarrow B$ be
the bijection that sends each block $C\in\tau$ to the unique
element $b$ in $C\cap B.$ For $(B,B')\in\pi\times\pi,$ we have
 \begin{equation}\label{condition}\Phi_{B,B'}^{\tau}= \Psi_{\tau,B'}\circ\Psi_{\tau,B}^{-1}.
 \end{equation}

Let $\Upsilon_U$ be the function
\begin{equation}\label{cloned1}
\Upsilon_U:(M\boxdot N)[U]\longrightarrow\mbox{ $M$-assemblies of
cloned $N$-structures on $U,$}
\end{equation}
that sends $(\pi,\tau,m,n)$ to $(\{n_B\}_{B\in\pi},\tau,m),$ where
$n_B=N[\Psi_{\tau,B}](n),$ for each $B\in\pi.$ From equation
(\ref{condition}) condition (\ref{cond}) is satisfied.

Let now $(\{n_B\}_{B\in\pi},\tau,m)$ be an $M$-assemblies of
cloned $N$-structures on $U$.
 From condition
(\ref{cond}) and equation (\ref{condition}), the $N$-structure
$n\in N[\tau],$ where $n:=N[\Psi_{\tau,B}^{-1}](n_B),$ remains the
same independently of the block $B\in \pi$ that we choose. It is
easy to check that $\Upsilon_U$ has as inverse the function that
sends $(\{n_B\}_{B\in\pi},\tau,m)$ to $(\pi,\tau,m,n).$

The family of bijections $\{\Upsilon_U\}_{U\in\mathbb{B}}$ is the
desired isomorphism.
\end{pf}

\subsection{Hyper-cloned rooted trees}

 The species of  $R$-enriched rooted trees could be
defined by the implicit combinatorial equation
\begin{equation}\label{root}
 \mathcal{A}_{R}=X\cdot R(\mathcal{A}_R).
\end{equation}
When $|R[\varnothing]|=1$ equation (\ref{root})
becomes
\begin{equation}\label{root2}
 \mathcal{A}_{R}=X+X\cdot R_+(\mathcal{A}_R).
\end{equation}
By changing the operation of substitution of species by the
arithmetic product in equation (\ref{root2}) we obtain the
combinatorial implicit equation for a new kind of structures, the
\emph{$R$-enriched hyper-cloned  rooted trees} ($R$-enriched
HRT's)
\begin{equation}\label{recursive}
 \mathcal{H}_{R}=X+X\cdot(R_+\boxdot \mathcal{H}_R).
\end{equation}
This equation leads to the following recursive definition: an
$\mathcal{H}_{R}$-structure on a set $U$ is either a singleton
vertex (when $|U|=1$), or is obtained by choosing a vertex in
$a_0\in U$(the root)  and attaching to it an $R_+$-assembly of
cloned $\mathcal{H}_R$-structures on $U\setminus \{a_0\}.$ To give
an explicit description of this kind of structures we need some
previous notation.

Let $t_U$ be an $R$-enriched rooted tree on $U$. The subset of $U$
formed by the non-leave elements will be denoted as $U^+.$ For
$a\in U^+$, the set $U_a$ will be the set of vertices in $U$ that
precede $a$ when the edges are oriented towards the root. The
 partition $\pi_a$ will be the partition of $U_a$ induced by the forest of
$R$-enriched trees that are attached to $a,$ and
$\{t_B^a\}_{B\in\pi_a}$ such forest.

\begin{Definition}
An \emph{$R$-enriched HRT} on a finite set $U$ is an $R$-enriched
tree $t_U$ together with a family of partitions, $\{\tau_a\}_{a\in
U^+}$, $\tau_a\in \mathrm{Par}[U_a]$, satisfying the conditions:
\begin{itemize}
\item[(i)]For every $a\in U^+,$ $(\pi_a,\tau_a)$ is a rectangle on
$U,$ \item[(ii)] For every pair $B, B'\in\pi_a,$
\begin{itemize}\item [(a)]$\Phi^{\tau_a}_{B,B'}:t_B^a\longrightarrow t_{B'}^a$ is an
isomorphism of $R$-enriched trees, \item[(b)]\label{homologia}If
$\Phi^{\tau_a}_{B,B'}(a_1)=a_2$ with $a_1\in B\cap U^+,$
 then $\Phi^{\tau_a}_{B,B'}(\tau_{a_1})=\tau_{a_2}.$
\end{itemize}
\end{itemize}
\end{Definition}

By condition (ii)(b), the family $\{\tau_a\}_{a\in U^+}$ is
completely determined by the partitions
$\{\tau_{a_0},\tau_{a_1},\dots,\tau_{a_k}\}$ on any branch
$\{a_0,a_1,\dots,a_k\}\subseteq U^+$ of $t_U.$
\begin{figure}[h]
  \centering
   \includegraphics[height=7.2cm]{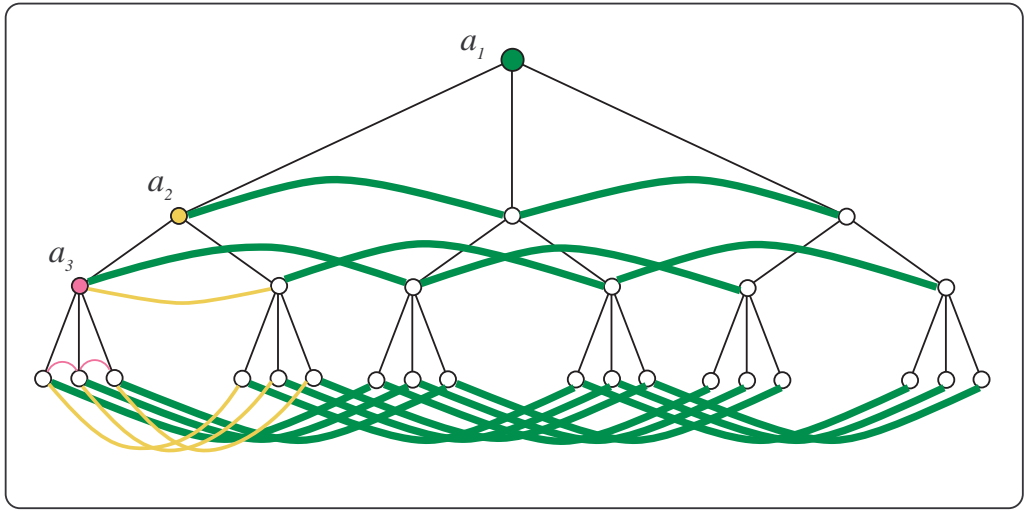}
   \caption{Element of $\mathcal{H}_{E}.$ Homologous elements according to
   $\tau_{a_1}, \tau_{a_2}, \tau_{a_3}$ are linked
    with different kinds of lines.}
   \label{hiper}
 \end{figure}

From equation (\ref{recursive}) we get the recursion:
\begin{eqnarray}
|\mathcal{H}_{R}[1]|&=&1,\\
|\mathcal{H}_{R}[n+1]|&=&\sum_{d|n}
\genfrac{\{}{\}}{0pt}{0}{n}{d}|R[d]||\mathcal{H}_{R}[n/d]|\; n\geq
1.
\end{eqnarray}

In particular, $\mathcal{H}_L$ is the species of rooted achiral
trees (see \cite{Ha2}). The ten first coefficients of the sequence
$|\mathcal{H}_E[n]|$ are shown in Table \ref{tab}.
\begin{table}[h!]
  \centering
\begin{tabular}{|c||c|c|c|c|c|c|c|c|c|c|}
  \hline
  $n$ & 1 & 2 & 3 & 4 & 5 & 6 & 7 & 8 & 9 & 10\\
  \hline
  $|\mathcal{H}_E[n]|$ & 1 & 1 & 2 & 3 & 10 & 11 & 192 & 193 & 3554 & 10080\\
  \hline
\end{tabular}
\caption{The ten first coefficients $|\mathcal{H}_E[n]|.$
}\label{tab}
\end{table}

\section{Multiplicative species}

We have  a notion of multiplicative species, the ``categorified''
analogous of the multiplicative arithmetic function in analytic
number theory (see \cite{Ap}).
\begin{Definition}
 Let $M$ be a species of structures satisfying the condition
 $M\left[\varnothing\right]=\varnothing.$ We say
 that $M$ is \emph{multiplicative} if,
 \begin{equation}
  M_{rs}=M_r\boxdot M_s,
 \end{equation}
 whenever $(r,s)=1.$
\end{Definition}

For example, using the isomorphisms of groups $\{1\}\simeq\{1\}\times\{1\}$ and
 $C_{rs}\simeq C_r\times C_s,$ when $(r,s)=1$, we have the species
 $L_+$ and $\mathcal{C}$ are multiplicative, respectively.
 The following three propositions are proved straightforwardly
 \begin{prop}\label{prop6}
 Let $M$ and $N$ be two multiplicative species of structures.
 Then the species $M^{\bullet}$ and
 $M\boxdot N$ are multiplicative.
\end{prop}

Combining the above examples and the above proposition, we find
that the species of regular octopuses $\mathcal{C}\boxdot L_+$ is
multiplicative.

\begin{prop}\label{mult_spec}
 Let $M$ be a species of structures. Then $M$ is multiplicative,
 if and only if, $M_1=X,$ and,
 for $n\geq 2,$
 \begin{equation}
  M_n=M_{p_1^{\alpha_1}}\boxdot M_{p_2^{\alpha_2}}\boxdot\cdots\boxdot M_{p_k^{\alpha_k}},
 \end{equation}
 with $n=\prod_{i=1}^k p_i^{\alpha_i}$ the canonical prime factorization of the
 integer $n.$
\end{prop}

\begin{prop}\emph{(}\textsc{Euler product formula}\emph{)}
 Let $M$ be a multiplicative species of structures. Then
 \begin{equation}
  M=\boxdot_{p\in \mathbf{P}}\left(X+M_p+M_{p^2}+\cdots\right),
 \end{equation}
 where $\mathbf{P}$ denotes the set of prime numbers.
\end{prop}

The following corollary follows immediately by taking generating
 series in the previous proposition.

 \begin{cor}
 Let $M$ be a multiplicative species of structures. Then:
 \begin{eqnarray}
 M(x)&=&\boxdot_{p\in
 \mathbf{P}}\left(x+M_p(x)+M_{p^2}(x)+\ldots\right),\\
\mathcal{D}_M(s)&=&\prod_{p\in\mathbf{P}}\left(1+\sum_{k\geq
1}\frac{|M[p^k]|}{p^k!p^{-ks}}\right),\\
  Z_M(x_1,x_2,\ldots)&=& \boxdot_{p\in \mathbf{P}}\left(x_1+Z_{M_p}(x_1,x_2,\ldots)+
      Z_{M_{p^2}}(x_1,x_2,\ldots)+\cdots\right).
 \end{eqnarray}
\end{cor}

\begin{Example}
For the multiplicative species $\mathcal{C}$ of cyclic
permutations we have the identities:
\begin{eqnarray}
 \mathcal{C}(x)&=&\boxdot_{p\in
 \mathbf{P}}\left(x+\sum_{k\geq 1}\frac{x^{p^k}}{p^k}\right)=\boxdot_{p\in
 \mathbf{P}}\left(x-\frac{x^p}{p}\right)^{\boxdot\langle -1\rangle},\\
\mathcal{D}_{\mathcal{C}}(s)&=&\zeta(s+1)=\prod_{p\in\mathbf{P}}\left(1-p^{-(s+1)}\right)^{-1},\\
  Z_{\mathcal{C}}(x_1,x_2,\ldots)&=& \boxdot_{p\in \mathbf{P}}\left((1-p^{-1})
  \sum_{j\geq 0}\sum_{k\geq
  j}\frac{x^{p^{k-j}}_{p^j}}{p^{k-j}}\right).
 \end{eqnarray}
\end{Example}

\section{Oligomorphic groups and species}
Let $A$ be an at most countable set, and $U$ a finite set. Denote
by $A^U$ and by $(A)_U$ the set of functions and injective
functions from $U$ to $A$ respectively. A permutation group $G$ on
the set $A$ is called {\em oligomorphic} if $G$ has only finitely
many orbits on $A^U$ for every finite set $U$.

 For $h\in A^U$ denote by $\ker(h)$ the partition of $U$ whose
 blocks are the non-empty pre-images of elements of $A$ by $h$.
 Recall that each function $f\in A^U$ can be identify with a pair
 $(\pi,\hat{h})$, where $\pi=\ker(h)$ and $\hat{h}\in (A)_{\pi}$
 is the injective function $\hat{h}(B):=h(b),$ for every $B\in \pi,
 $ $b$ being an arbitrary element of $B$.

\begin{Definition} Let $G$, $A$ and $U$ be as above.
Define the species of structures $F^*_G$ by $F^*_G[U]=A^U/G$, the
(finite) set of orbits of $A^U$ under the action of $G$. For a
bijection $\sigma:U\longrightarrow V,$ define the bijection
 \begin{equation}\label{bijolig}
  \begin{array}{cl}
    F^*_G[\sigma] : &\!\!\! F^*_G[U] \longrightarrow F^*_G[V]\\
                    &\,\,\,\,\,\,\,\,\, \overline{f}\, \longmapsto\,
                    \overline{f\circ\sigma^{-1}}.
  \end{array}
 \end{equation}
 This bijection is well defined since $\sigma$ commutes with the action of
 $G$ over $A^U.$ In an analogous way we define the species $F_G$
 of $G$-orbits of injective functions.
\end{Definition}
Recall that for a finite set $A$ and a species $M,$ the species of
$M$-enriched functions is denoted by $M^A$ (see \cite{Jo1}).
Observe that when $A$ is finite and $G$ is the identity subgroup
of $S_A,$ $F_G[U]=A_U$ is isomorphic to $(1+X)^A$ and
$F_G^*[U]=A^U$ is
isomorphic to $E^A.$\\

\noindent \textbf{Remark.}\, Cameron \cite{Ca2} has studied the
three following counting problems : how many elements in (a)
$F_G[n],$ (b) $\widetilde{F_G}[n],$ (c) $F_G^*[n]?$ Equivalently,
how many $G$-orbits in (a) $n$-tuples of distinct elements, (b)
$n$-sets, (c) all $n$-tuples?

\begin{prop}
We have the following combinatorial identity
\begin{equation}
F_G^*=F_G(E_+).
\end{equation}
\end{prop}
\begin{pf} Let $h\in A^U$, since the action of $G$ does not
affect the kernel of $h$, the orbit $\overline{h}$ of $h$ can be
identified with the pair $(\pi,\overline{\hat{h}})$, where
$\pi=\ker(h).$ Obviously $\overline{\hat{h}}\in F_G[\pi]$. This
defines a natural bijection from $F_G^*[U]$ to $F_G(E_+)[U].$
\end{pf}

\subsection{The modified arithmetic product}
Let $G$ and $H$ be two oligomorphic groups of permutations on the
sets $A$ and $B$ respectively. In \cite{Ca} Cameron et al.\! deal
with the enumerative problem in the remark above for the product
group $G\times H$ acting over $A\times B$. We now introduce the
analogous problem in the more general context of species of
structures.

\begin{Definition}(\textsc{Modified arithmetic product of species})
Let $M$ and $N$ be two species of structures. Denote by
$P_{\mathcal{R}}$ the species of partial rectangles. We define the
modified arithmetic
 product of $M$ and $N$ by
\begin{equation}
(M\widetilde{\boxdot}N)[U]:=\sum_{(\pi,\tau)\in
P_{\mathcal{R}}[U]}M[\pi]\times N[\tau],
\end{equation}
where the sum represents the disjoint union and $U$ is a finite
set. For a bijection $\sigma:U\longrightarrow V,$ the transport
$(M\widetilde{\boxdot} N)[\sigma]$ is as in Definition \ref{def2}.
\end{Definition}
\begin{figure}[ht]
  \centering
   \includegraphics[height=8cm]{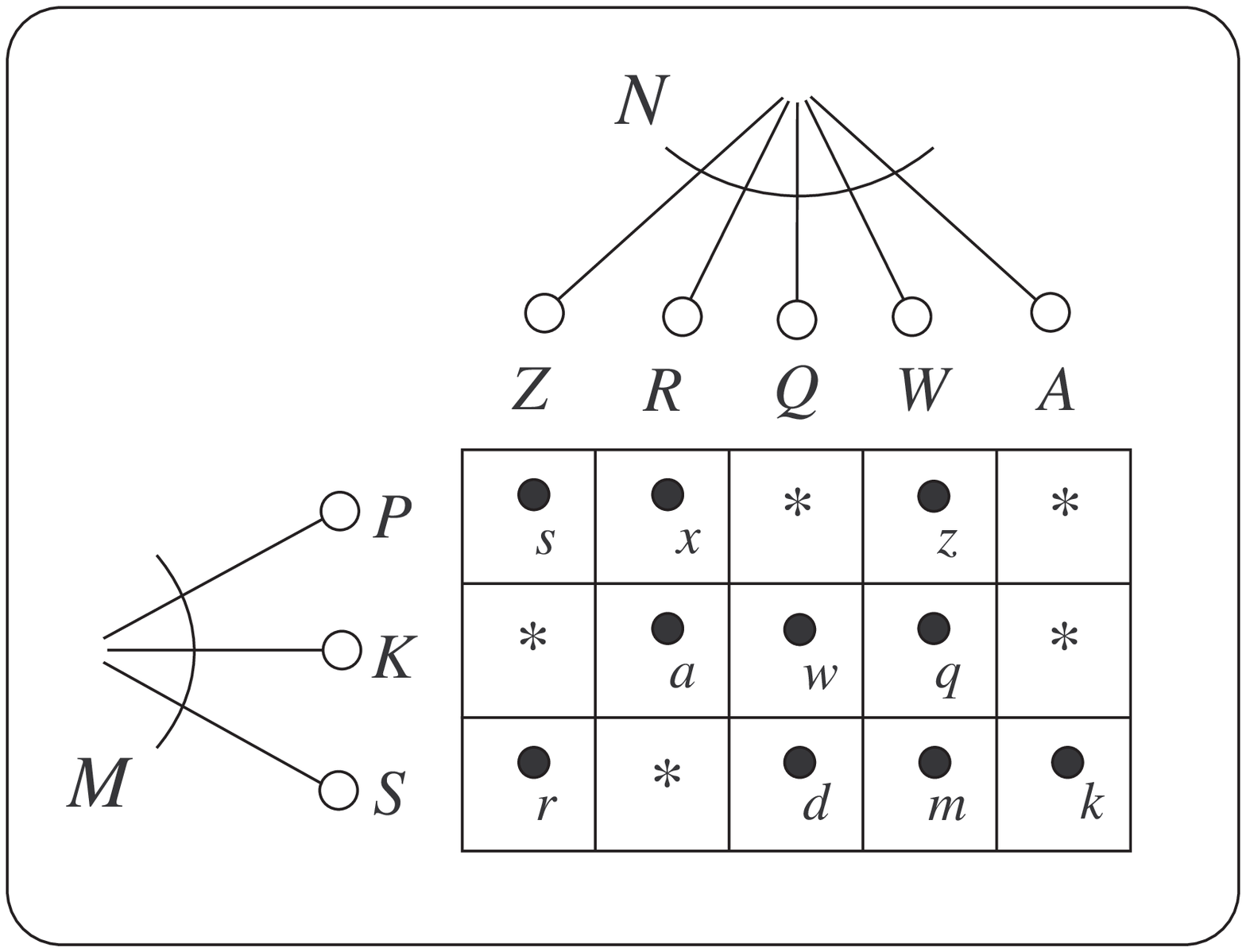}
   \caption{An $(M\widetilde{\boxdot} N)$-structure on a 10-set.}
   \label{modproduct}
 \end{figure}
Some of the properties of the arithmetic product in Proposition
\ref{prop} have their analogous in this context.
\begin{prop}\label{prop1}
 Let $M,N$ and $R$ be species of structures.
 The product $\widetilde{\boxdot}$ has the following properties:
 \begin{eqnarray}
    M\widetilde{\boxdot} N&=&N\widetilde{\boxdot} M,\\
    \label{asociatividad1}M\widetilde{\boxdot}(N\widetilde{\boxdot} R)
    &=&(M\widetilde{\boxdot} N)\widetilde{\boxdot} R,\\
    M\widetilde{\boxdot}(N+R)&=&M\widetilde{\boxdot} N+M\widetilde{\boxdot} R,\\
    X\widetilde{\boxdot}M&=&M_+,\\
    M\widetilde{\boxdot} (1+X)&=&(X+1)\widetilde{\boxdot} M=M,\\
    \label{principal}M\widetilde{\boxdot}(1+X)^{[n]}&=&M((1+X)^{[n]}_+).
 \end{eqnarray}
\end{prop}
\begin{pf}  We will only prove identity (\ref{principal}). An
element of $(M\widetilde{\boxdot}(1+X)^{[n]})[U]$ is of the form
$(\pi,\tau,m,\hat{f}),$ where $(\pi,\tau)$ is a partial rectangle
on $U,$ $m\in M[\pi],$ and $\hat{f}:\tau\rightarrow [n]$ is an
injective function. Recall that the pair $(\tau,\hat{f})$ can be
identify with a function $f:U\rightarrow [n]$ whose kernel is
$\tau$. Since $\pi\wedge\tau=\hat{0},$ the restriction $f_B$ of
$f$ to each block $B$ of $\pi$ is injective. Conversely, if all
the functions in a family $\{f_B\}_{B\in\pi}$ are injective, then
$\pi\wedge\tau=\hat{0}$, $\tau$ being the kernel of
$f:=\cup_{B\in\pi}f_B.$ Then, the correspondence
\begin{eqnarray*}
\Omega_U:(M\widetilde{\boxdot}(1+X)^{[n]})[U]&\longrightarrow &
M((1+X)^{[n]}_+)[U]\\
(\pi,\tau,m,\hat{f})&\longmapsto & (\{f_B\}_{B\in\pi},m),
\end{eqnarray*}
is a natural bijection.
\end{pf}

Like in equation (\ref{general}) the product of a family $\{ M_i\}_{i=1}^k$ of
species of structures is given by
\begin{equation}
(\widetilde{\boxdot}_{i=1}^k
M_i)[U]=\sum_{(\pi_1,\pi_2,\dots,\pi_k)\in
P_{\mathcal{R}}^{(k)}[U]}\prod_{i=1}^k M_i[\pi_i],
\end{equation}
where $P_{\mathcal{R}}^{(k)}$ is the species of $k$-partial
rectangles.

We have the following Theorem.
\begin{theo}\label{product} Let $G$ and $H$ be two oligomorphic
groups acting on sets $A$ and $B$ respectively. Then
\begin{equation}
F_{G\times H}=F_G\widetilde{\boxdot}F_H.
\end{equation}
\end{theo}
\begin{pf}
Let $h:U\longrightarrow A\times B$ be an injective function. Let
$h_1:U\longrightarrow A$ and $h_2:U\longrightarrow B$ be its
components, i.e. $h(u)=(h_1(u),h_2(u))$ for $u\in U.$ Let
$\pi=\ker(h_1)$ and $\tau=\ker(h_2)$. It is clear that
$\ker(h)=\ker(h_1)\wedge \ker(h_2)=\pi\wedge\tau$, and since $h$
is injective, $\pi\wedge\tau=\hat{0}.$ Then $h$ can be identify
with the tuple $(\pi,\tau,\hat{h}_1,\hat{h}_2),$ and its orbit
under the action of $G\times H$ with
$(\pi,\tau,\overline{\hat{h}}_1,\overline{\hat{h}}_2),$ where
$\overline{\hat{h}}_1\in F_G[\pi]$ and $\overline{\hat{h}}_2\in
F_H[\tau]$. This defines a natural bijection between $F_{G\times
H}[U]$ and $(F_G\widetilde{\boxdot}F_H)[U].$
\end{pf}

Take $A=[m],B=[n],$ and $G,H$ being the identity subgroups of
$S_m$ and $S_n$ respectively. We obtain the isomorphism
\begin{equation}\label{segunda}
(1+X)^{[m]}\widetilde{\boxdot}(1+X)^{[n]}=(1+X)^{[n]\times[m]}.
\end{equation}

 The exponential generating series of the
modified arithmetic product of species of structures is not as straightforward
to compute as in the arithmetic product case. However, the
identity
\begin{equation}\label{Cameron}
F_{G\times H}^*(x)=F_G^*(x)\times F_H^*(x),
\end{equation}
proved in \cite{Ca}, provides a device to compute this series,
\begin{equation}
 (F_G\widetilde{\boxdot}F_H)(x)=F_{G\times H}(x).
\end{equation}
It has motivated the following general combinatorial
identity.
\begin{theo}
Let $\{M_i\}_{i=1}^k$ be a family of species of structures. Then we have
\begin{equation}\label{canf}
\widetilde{\boxdot}_{i=1}^kM_i(E_+)=\times_{i=1}^kM_i(E_+),
\end{equation}
where $\times$ is the operation of cartesian product of species,
\begin{equation}
(M\times N)[U]=M[U]\times N[U].
\end{equation}
\end{theo}
\begin{pf}
It is enough to prove the identity for $k=2.$ Consider the species
$\mathrm{Par}$ of set partitions. For a finite set $U$, let
$\leq_U$ be the refinement order on $\mathrm{Par}[U]$. For
$\eta\in \mathrm{Par}[U]$, let $\mathrm{C}_{\eta}$ denote the
order coideal of $(\mathrm{Par}[U],\leq_U)$ of the elements
greater than or equal to $\eta.$ For $\pi\in\mathrm{C}_{\eta}$ and
$B\in\pi,$ let $\hat{B}=\{C\in\eta\mid C\subseteq B\}$ and
$\hat{\pi}=\{\hat{B}\mid B\in\pi\}.$ Clearly, $\hat{\pi}$ is a
partition of $\eta$, and it is easy to see that the correspondence
\begin{eqnarray*}
\mathrm{C}_{\eta}&\longrightarrow&
(\mathrm{Par}[\eta],\leq_{\eta})\\\pi&\longmapsto&\hat{\pi}
\end{eqnarray*}
is an order isomorphism. Then the partition $\pi\wedge\tau=\eta$
is an element of $\mathrm{C}_{\eta}$ if and only if
$(\hat{\pi},\hat{\tau})$ is a partial rectangle on $\eta.$ The
right hand side of (\ref{canf}), for $k=2$, evaluated in a set $U$
is equal to
\begin{eqnarray}
(M_1(E_+)\times M_2(E_+))[U]&=&\left(\sum_{\pi\in
\mathrm{Par}[U]}M_1[\pi]\right)\times \left(\sum_{\tau\in
\mathrm{Par}[U]}M_2[\tau]\right)\\&=&\sum_{(\pi,\tau)\in
\mathrm{Par}[U]\times \mathrm{Par}[U]}M_1[\pi]\times
M_2[\tau]\\\label{ultimo}&=&\sum_{\eta\in
\mathrm{Par}[U]}\sum_{\pi\wedge\tau=\eta}M_1[\pi]\times M_2[\tau].
\end{eqnarray}
For any partition $\varrho$ in $\mathrm{C}_{\eta},$ let
$f_{\varrho,\eta}:\varrho\longrightarrow \hat{\varrho}$ be the
bijection sending each block $B$ of $\varrho$ to $\hat{B}$. The
family of bijections
\begin{equation}\alpha_U:\sum_{\eta\in
\mathrm{Par}[U]}\sum_{\pi\wedge\tau=\eta}M_1[\pi]\times
M_2[\tau]\longrightarrow\sum_{\eta\in
\mathrm{Par}[U]}\sum_{(\hat{\pi},\hat{\tau})\in
P_{\mathcal{R}}[\eta]}M_1[\hat{\pi}]\times M_2[\hat{\tau}]
\end{equation}
\begin{equation}
\alpha_U:=\sum_{\eta\in
\mathrm{Par}[U]}\sum_{\pi\wedge\tau=\eta}M_1[f_{\pi,\eta}]\times
M_2[f_{\tau,\eta}]\end{equation}defines a natural transformation
\begin{equation}
\alpha:M_1(E_+)\times M_2(E_+)\longrightarrow
(M_1\widetilde{\boxdot}M_2)(E_+).
\end{equation}
\end{pf}

Taking exponential generating series and cycle index series in
identity (\ref{canf}), we obtain the following
\begin{cor}Let $\{M_i\}_{i=1}^k$ be as above. Then the following
generating function identities hold
\begin{eqnarray}
\label{clumsy1}(\widetilde{\boxdot}_{i=1}^kM_i)(e^x-1)
&=&\times_{i=1}^k M_i(e^x-1),\\
Z_{\widetilde{\boxdot}_{i=1}^kM_i}(\mathbf{x})\ast
Z_{E_+}(\mathbf{x}) &=&\times_{i=1}^k(Z_{M_i}(\mathbf{x})\ast
Z_{E_+}(\mathbf{x})),
\end{eqnarray}
where $\times$ means coefficient-wise product or Hadamard product
as in \emph{\cite{BLL}}, $\ast$ means plethystic substitution, and
\begin{equation}
Z_{E_+}(\mathbf{x})=\exp\left(\sum_{n\geq
1}\frac{x_n}{n}\right)-1.
\end{equation}
\end{cor}

Using equation (\ref{clumsy1}) with $M_i=E$, for $i=1,\dots, k,$
we recover the first identity of Theorem 1 in \cite{Can},
\begin{equation}\label{B}
 P_\mathcal{R}^{(k)}(e^x-1)=(e^{e^x-1})^{\times
k}=\sum_{n\geq 0}(B_n)^k\frac{x^n}{n!},
\end{equation}
where $B_n$ is the $n$-th Bell number, the number of partitions of
the set $[n].$

 In order to have the identity
$(M\widetilde{\boxdot}N)(x)=M(x)\widetilde{\boxdot}N(x)$ for
 two species of structures
$M$ and $N$, following (\ref{clumsy1}) we make the following
Definition.
\begin{Definition}
For two formal power series $F(x)$ and $G(x)$ define the product
$\widetilde{\boxdot}$ by
\begin{equation}\label{clumsy2}
F(x)\widetilde{\boxdot}G(x)=\left(F(e^x-1)\times
G(e^x-1)\right)\circ(\ln(1+x)).
\end{equation}
\end{Definition}
It is easy to see that this product is commutative and
distributive with respect to the sum,
\begin{equation}
\left(F(x)+H(x)\right)\widetilde{\boxdot}G(x)
=F(x)\widetilde{\boxdot}G(x)+H(x)\widetilde{\boxdot}G(x).
\end{equation}

\subsection{The shift trick}
Sometimes the equation (\ref{clumsy2}) is too clumsy to make
computations. We will provide a more efficient method. Previous to
that we need the following Lemma.
\begin{Lemma}Let $F(x)$ be a formal power series.
For $m$ and $n$ nonnegative integers we have the
identities:
\begin{eqnarray}\label{trick}\label{trick2}F(x)\widetilde{\boxdot}(1+x)^n&=&F((1+x)^n-1),\\
\label{trick1}(1+x)^m\widetilde{\boxdot}(1+x)^n&=&(1+x)^{mn}.
\end{eqnarray}
\end{Lemma}
\begin{pf}
Equation (\ref{trick2}) follows from identity (\ref{principal}).
Equation(\ref{trick1}) follows from (\ref{trick2}) or by taking
generating functions in (\ref{segunda}).
\end{pf}

From this lemma we recover the following result of Pittel
\cite{Pi}.
\begin{prop}Let $k$ be a fixed positive integer. The exponential
generating series $P_{\mathcal{R},k}(x)$, of the number of
partial rectangles $(\pi,\tau)$, with $|\pi|=k$, is
\begin{equation}
P_{\mathcal{R},k}(x)=\frac{1}{k!e}\sum_{l\geq
0}\frac{1}{l!}\left((x+1)^l-1\right)^k.
\end{equation}
\end{prop}
\begin{pf}
The required species is $P_{\mathcal{R},k}:=E_k\widetilde{\boxdot} E.$
 Its exponential generating series is
\begin{equation}
(E_k\widetilde{\boxdot}
E)(x)=\frac{x^k}{k!}\widetilde{\boxdot}e^x=\frac{1}{k!e}x^k\widetilde{\boxdot}e^{(x+1)}=\frac{1}{k!e}
\sum_{l\geq 0}\frac{1}{l!}x^{k}\widetilde{\boxdot}(x+1)^l.
\end{equation}
Use equation (\ref{trick2}) to finish the proof.
\end{pf}

 The algorithm to compute the product
$F(x)\widetilde{\boxdot}G(x),$ of two generating series $F(x)$ and
$G(x),$ runs as follows:
\begin{enumerate}
\item Express $F(x)=F_1(x+1)$ and $G(x)=G_1(x+1)$ as power series
of $(x+1),$
\item use the distributive property and
equation (\ref{trick1}) to compute
$$H(x+1)=F_1(x+1)\widetilde{\boxdot}G_1(x+1),$$
\item express back
$H(x+1)$ as a power series of $x.$
\end{enumerate}

 Now we solve some enumerative problems.
\begin{theo}
The number $M(m,n,r)$of $m\times n$ $(0,1)$-matrices with exactly
 $r$ entries equal to $1$ and no zero row or columns, is given by
\begin{equation}
M(m,n,r)=\sum_{l\geq r}\sum_{d|l}(-1)^{n+m-(d+l/d)}\genfrac{(}{)
}{0pt}{0}{m}{d}\genfrac{(}{)}{0pt}{0}{n}{l/d}
\genfrac{(}{)}{0pt}{0}{l}{r}.
\end{equation}
\end{theo}
\begin{pf}
The structures of the species $X^n\widetilde{\boxdot}X^m$ are the
linearly ordered $m\times n$ partial rectangles. A structure of
$(X^n\widetilde{\boxdot}X^m)[r]$ can be thought of as a $m\times n$
matrix with entries $1,2,\dots,r$, without repetitions, zero
elsewhere, and no zero row or columns. Then, $M(m,n,r)$ is the
coefficient of $x^r$ in the generating series
\begin{equation}
(X^m\widetilde{\boxdot}X^n)(x)=x^m\widetilde{\boxdot}x^n.
\end{equation}
By shifting we get
\begin{eqnarray}
x^m\widetilde{\boxdot}x^n&=&(x+1-1)^m\widetilde{\boxdot}(x+1-1)^n\\
&=&\sum_{j,k}\genfrac{(}{) }{0pt}{0}{m}{j}\genfrac{(}{)
}{0pt}{0}{n}{k}(-1)^{m+n-(j+k)}(x+1)^j\widetilde{\boxdot}(x+1)^k\\
&=&\sum_{j,k}\genfrac{(}{) }{0pt}{0}{m}{j}\genfrac{(}{)
}{0pt}{0}{n}{k}(-1)^{m+n-(j+k)}(x+1)^{jk}\\
&=&\sum_{j,k,r}\genfrac{(}{) }{0pt}{0}{m}{j}\genfrac{(}{)
}{0pt}{0}{n}{k}\genfrac{(}{) }{0pt}{0}{jk}{r}(-1)^{m+n-(j+k)}x^r.
\end{eqnarray}
Making the change $l=jk,$ we obtain the result.
\end{pf}

\begin{cor}The number $|P_{\mathcal{R},m,n}[r]|$ of $m\times n$ partial rectangles
on $r$ elements, is given by
\begin{equation} |P_{\mathcal{R},m,n}[r]|=\sum_{l\geq
r}\sum_{d|l}\genfrac{\{}{\}
}{0pt}{0}{l}{d}\frac{(-1)^{m+n-(d+l/d)}}{(m-d)!(n-l/d)!(l-r)!}.
\end{equation}
\end{cor}
\begin{pf}
$|P_{{\mathcal R},m,n}[r]|$ is equal to
$|(E_m\widetilde{\boxdot}E_n)[r]|,$ which is the coefficient of
$\displaystyle\frac{x^r}{r!}$ in the generating series
\begin{equation}
(E_m\widetilde{\boxdot}E_n)(x)=\frac{x^m}{m!}\widetilde{\boxdot}\frac{x^n}{n!}.
\end{equation}
Then,
$$|P_{\mathcal{R},m,n}[r]|=r!\frac{M(m,n,r)}{m!n!}=\sum_{l\geq
r}\sum_{d|l}\genfrac{\{}{\}
}{0pt}{0}{l}{d}\frac{(-1)^{m+n-(d+l/d)}}{(m-d)!(n-l/d)!(l-r)!}.$$
\end{pf}

We now give a very short direct proof of the beautiful formula
obtained by Pittel \cite{Pi}.
\begin{theo}\label{beauty}
The number $|P_{\mathcal{R}}^{(k)}[n]|,$ of\, $k$-tuples of
partitions $(\pi_1,\pi_2,\dots, \pi_k)$ on $[n]$ satisfying
$\pi_1\wedge\pi_2\wedge\dots\wedge\pi_k=\hat{0},$ is given by
\begin{equation}\label{kpartial}
|P_{\mathcal{R}} ^{(k)}[n]|=e^{-k}\sum_{i_1,i_2,\dots,i_k\geq
1}\frac{(i_1\cdots i_k)_n}{i_1!\cdots i_k!},
\end{equation}
where $(m)_n=m(m-1)\cdots(m-n+1).$
\end{theo}
\begin{pf}
 By the definition of $\widetilde{\boxdot}$-product,
$$
P_{\mathcal{R}}^{(k)}=E^{\widetilde{\boxdot}
k}=\underbrace{E\widetilde{\boxdot}\cdots\widetilde{\boxdot}E}_{k
\rm\; factors}.
$$
The exponential generating series of this species is
$(e^x)^{\widetilde{\boxdot}k}.$ Following the algorithm, we have
$e^x=e^{-1}e^{(x+1)}.$ Then
\begin{eqnarray}
E^{\widetilde{\boxdot}
k}(x)&=&\left(e^{-1}e^{(x+1)}\right)^{\widetilde{\boxdot}k}\\
     &=&e^{-k}\sum_{i_1,i_2,\dots,i_k\geq 0}
\frac{(x+1)^{i_1}}{i_1!}\widetilde{\boxdot}\cdots
\widetilde{\boxdot}
 \frac{(x+1)^{i_k}}{i_k!}\\
 &=&
 e^{-k}\sum_{i_1,i_2,\dots,i_k\geq 0}\frac{(x+1)^{i_1i_2\cdots i_k}}{i_1!i_2!\cdots
 i_k!}\\
 &=&
 \sum_{n\geq 0}\left(e^{-k}\sum_{i_1,i_2,\dots,i_k\geq 0}\frac{(i_1\cdots
 i_k)_n}{i_1!i_2!\cdots i_k!}\right)\frac{x^n}{n!}.
\end{eqnarray}
\end{pf}

As a corollary, we obtain an remarkable identity
\begin{cor}For $n\geq 1,$
\begin{equation}
|P_{\mathcal{R}}^{(k)}[n]|=e^{-k}\sum_{l\geq
n}\frac{|\mathcal{R}^{(k)}[l]|}{(l-n)!}.
\end{equation}
\end{cor}
\begin{pf}
Making the change $l=i_1i_2\cdots i_k$ in equation (\ref{kpartial})
we obtain
\begin{eqnarray}
|P_{\mathcal{R}}^{(k)}[n]| &=&e^{-k}\sum_{l\geq n}\sum_{i_1
i_2\cdots i_k\geq l}
\frac{(l)_n}{i_1!\cdots i_k!}\\
&=&e^{-k}\sum_{l\geq n}\frac{1}{(l-n)!} \sum_{i_1 i_2\cdots
i_k\geq l}\frac{l!}{i_1!\cdots i_k!}.
\end{eqnarray}
To finish the proof we recall equation (\ref{krectangulo}).
\end{pf}

Theorem \ref{beauty} is a particular case of the following general
result, that can be proved without much extra effort.

\begin{theo}Let $\{M_i\}_{i=1}^k$ be a family of species of structures whose
exponential generating series, $M_i(x)=F_i(x+1)$, expressed as
power series of $(x+1),$ are given,
\begin{equation}
M_i(x)=F_i(x+1)=\sum_{n\geq 0} b_{n}^{(i)}\frac{(x+1)^n}{n!},\quad
i=1,\dots,k.
\end{equation}
Then,
\begin{eqnarray}
(\widetilde{\boxdot}_{i=1}^k M_i)(x)&=&\sum_{n\geq
0}\left(\sum_{i_1,i_2,\dots,i_k\geq 0} b_{i_1}^{(1)}\cdots
b_{i_k}^{(k)}\frac{(i_1\cdots
 i_k)_n}{i_1!i_2!\cdots
 i_k!}\right)\frac{x^n}{n!}\\
 &=&a_0+\sum_{n\geq 1}\left(\sum_{l\geq
n}\frac{1}{(l-n)!}\left[\frac{x^n}{n!}\right]
 (F_1\boxdot\cdots\boxdot F_k)(x)\right)\frac{x^n}{n!},
\end{eqnarray}where
$a_0=|(\widetilde{\boxdot}_{i=1}^k M_i)[0]|
=\prod_{i=1}^k|M_i[0]|.$
\end{theo}

\begin{Example}
Given the expansions:
\begin{eqnarray}
L(x)&=&\frac{1}{2\left(1-\frac{x+1}{2}\right)}=\sum_{n\geq 0}\frac{(x+1)^n}{2^{n+1}},\\
\mathcal{C}(x)&=&\ln\left(\frac{1}{2\left(1-\frac{x+1}{2}\right)}\right)
=\ln(2^{-1})+\sum_{n\geq
1} \frac{(x+1)^n}{n2^{n}}.
\end{eqnarray}
We obtain formulas for:
\begin{itemize}
\item The number of\, $(0,1)$-matrices
 of any size, with $n$ ones, and
with no zero row or column,
\begin{equation}
\frac{|(L{\widetilde{\boxdot}}L)[n]|}{n!}=\frac{1}{4n!}\sum_{r,s\geq
0}\frac{(rs)_n}{2^{r+s}}.
\end{equation}
\item The number of matrices of any size up to column
permutations, with $n$ different elements, zero elsewhere and with
no zero row or column,
\begin{equation}
|(L\widetilde{\boxdot}E)[n]|=\frac{1}{2e}\sum_{r,s\geq
0}\frac{(rs)_n}{2^{r}s!}.
\end{equation}
\item The number of linearly ordered $k$-partial rectangles on
$[n]$,
\begin{equation}
|L^{\widetilde{\boxdot}k}[n]|=\frac{1}{2^k}\sum_{i_1\cdots i_k\geq
0}\frac{(i_1\cdots i_k)_n}{2^{i_1+\cdots+i_k}}.
\end{equation}
\item The number of cyclic $k$-partial rectangles on $[n]$,
\begin{equation}
|\mathcal{C}^{\widetilde{\boxdot}k}[n]|=\sum_{i_1\cdots i_k\geq
1}\frac{(i_1\cdots i_k-1)_{n-1}}{2^{i_1+\cdots+i_k}}.
\end{equation}
\end{itemize}
\end{Example}

 \vspace{2cm}

\begin{tabular}[t]{c@{\extracolsep{5em}}c}
\textsc{Manuel Maia}                & \textsc{Miguel M\'endez}  \\
   Departamento de Matem\'atica  &  Departamento de
Matem\'atica
 \\
   Facultad de Ciencias          &  IVIC  \\
   Universidad Central de Venezuela &  Carretera Panamericana, Km 11 \\
   Av. Los ilustres, Los Chaguaramos              &
Altos de Pipe, Estado Miranda\\
   A.P.: 20513, Caracas 1020--A             &
A.P. 21827,  Caracas 1020--A\\
   Venezuela            &
Venezuela\\
   \texttt{mmaia@euler.ciens.ucv.ve}         &
and\\
&
Departamento de Matem\'atica\\
           &
Universidad Central de Venezuela\\
          &
\texttt{mmendez@cauchy.ivic.ve}\\
\end{tabular}}

\end{document}